\documentclass{article}
\usepackage{amssymb}

\usepackage{graphicx} 
\usepackage[english]{babel}
\usepackage{amsthm}
\usepackage{tikz}

\newtheorem{definition}{Definition}
\newtheorem{manthminner}{Theorem}
\newenvironment{manthm}[1]{%
  \IfBlankTF{#1}
    {\renewcommand{\themanthminner}{\unskip}}
    {\renewcommand\themanthminner{#1}}%
  \manthminner
}{\endmanthminner}

\newtheorem{mandefinner}{Definition}
\newenvironment{mandef}[1]{%
  \IfBlankTF{#1}
    {\renewcommand{\themandefinner}{\unskip}}
    {\renewcommand\themandefinner{#1}}%
  \mandefinner
}{\endmandefinner}

\newtheorem{manleminner}{Lemma}
\newenvironment{manlem}[1]{%
  \IfBlankTF{#1}
    {\renewcommand{\themanleminner}{\unskip}}
    {\renewcommand\themanleminner{#1}}%
  \manleminner
}{\endmanleminner}
\usepackage{amsmath}
\usepackage{sectsty}
    \sectionfont{\large\centering}

\title{Degree 2 vertices in minimal prime graph complements}
\author{
  Bryan Alvarez, Micah Dorton, Thomas Michael Keller,\\ Lawrence Liu, Evan Zhang
}
\date{}

\usetikzlibrary{shapes.geometric, positioning}

\begin{document}

\maketitle
\begin{abstract}
Minimal prime graphs are connected graphs on at least two vertices whose complements satisfy the following 
conditions: triangle-freeness, 3-colorability, and edge-maximality with respect to the latter two properties. 
These graphs are prime graphs (or Gruenberg-Kegel graphs) of finite solvable groups with the maximum number of Frobenius actions among their
Sylow subgroups, and as such minimal prime graph complements have been shown to be highly structured, including, for instance, 
the presence of induced 5-cycles. It is also known that the minimum degree of minimal prime graph complements is 2. In this note, we show
that the existence of a degree 2 vertex in a minimal prime graph complement determines its whole structure: it is simply a 5-cycle with three 
vertices, exactly two of which are adjacent to each other, being duplicated finitely often. In particular, such graphs belong to a class of graphs known as reseminant.

\end{abstract}

\section*{Introduction}

The prime graph of a finite group \( G \), also known as its Gruenberg–Kegel graph, is defined on the set of prime divisors 
of \(|G|\), with edges connecting primes \( p \) and \( q \) if and only if \( pq \) divides the order of some element in \( G \). 
Introduced by Gruber et al.~\cite{gruber}, minimal prime graphs are prime graphs of solvable group with a 
maximal number of Frobenius actions among their Sylow subgroups. Technically spoken, this amounts to the following definition:
\begin{definition}
    A minimal prime graph (MPG) $\Gamma$ is a graph with $2$ or more vertices satisfying the following properties:
\begin{enumerate}
    \item $\Gamma$ is connected
    \item $\overline{\Gamma}$ is triangle-free
    \item $\overline{\Gamma}$ is $3$-colorable
    \item Removing any edge from $\Gamma$ will violate $(2)$ or $(3)$. This is the condition for minimality, which is equivalent to maximality in the complement.
\end{enumerate}
\end{definition}

It was already noticed in ~\cite{gruber} that vertex duplication, that is adding a vertex to a minimal prime graph by connecting it to one vertex and all its neighbors, yields a new minimal prime graph. It was also noted there that
minimal prime graph complements have minimal degree at least 2 and that
they always contain an induced 5-cycle. Florez et al.~\cite{florez} explored their diameters, Hamiltonian cycles, and self-complementarity, while Huang et al.~\cite{huang} introduced new methods for generating minimal prime graphs, including constructions beyond vertex duplication. These developments have led to the definition of \emph{reseminant graphs}, which are minimal prime graphs obtained by repeated vertex duplication from the 5-cycle \( C_5 \), the smallest minimal prime graph.
\\
In this paper, inspired by a question in \cite{joshua}, we prove that if a minimal prime graph complement contains a vertex of degree 2, then the graph is reseminant. In fact, it is even more special: it is a 5-cycle with three
vertices, exactly two of which are adjacent to each other, being duplicated finitely often. This shows that, perhaps surprisingly,  degree 2 vertices only show up in very restricted situations. \\
Before we dive into the proof, we would like to point out an interesting result by Goddard and Kleitman~\cite{goddard} dealing with a similar situation than we do in this paper: They study \emph{maximal triangle-free graphs}—graphs in which the addition of any edge creates a triangle. Their work establishes that such graphs contain large similarity classes (sets of vertices with identical neighborhoods). More precisely, they give a lower bound for the size of the largest similarity class in terms of the minimum degree and the
degrees of the vertices in the similarity class. They then use this to give a simplified proof of a conjecture by Woodall: that a graph with binding number at least \( \frac{3}{2} \) must contain a triangle.
\\ 
The connection to minimal prime graphs is obvious: the complement of a minimal prime graph is 3-colorable, triangle-free and edge-maximal with respect to these two conditions, thus "close" to a maximal triangle-free graph. Moreover, the notion of large similarity classes in maximal triangle-free graphs is closely related to vertex duplication in that a vertex and all its duplicates form a similarity class. Our main result implies that if in a minimal prime graph complement the minimal degree is 2, then there exist large similarity class
(of roughly a third of the vertices). This, however, does not follow from the Goddard-Kleitman result.

\section*{Results and proofs}
We first formally state the main result to be proved:

\begin{manthm} {2} \label{main}
    If there exists a vertex of degree 2 in a minimal prime graph complement $G$, then there exists an induced 5-cycle $S$ in $G$ and three vertices in $S$, exactly two of which are adjacent to each other, such that $G$ is obtained from $S$ by duplicating each of the three verices a suitable number of times. In particular, $G$ is reseminant.
\end{manthm}

Note that vertex duplication in the complement amounts to
choosing a vertex $v$, adding a new vertex $w$ and connecting $w$ 
to all the neighbors of $v$ (but not to $v$ itself). With this it
is easy to see that complements of reseminant prime graphs with minimal degree 2 can only be of the form described in Theorem \ref{main}.\\

Before we proceed, we explain our notation:
\begin{itemize}
    \item $V_G$ and $E_G$ are defined to be the vertex set and edge set of the graph $G$, respectively.
    \item An undirected edge between vertices $u$ and $v$ is denoted by $(u,v)$.
    \item $\overline{\Gamma}$ denotes the complement of the graph $\Gamma$.
    \item The degree of a vertex $v$ is denoted by $deg(v)$.
\end{itemize}

We now embark on proving Theorem \ref{main} using a number of intermediate
lemmas.\\

The theorem is of interest because 2 is the smallest possible degree a vertex of a minimal prime graph complement can have. This was proved as a consequence of a more general result in \cite[Corollary 11]{joshua}, but because of its importance and 
for the sake of completeness, we have decided to include here a direct proof of this fact. 

\begin{manlem}{3} \label{noDeg1}
    In a minimal prime graph complement, there cannot exist a vertex of degree one.
\end{manlem}

\begin{proof}
    Assume for the sake of contradiction that there exists minimal prime graph complement $G$, there exists $v\in V_G$ s.t. $deg(v)=1$. Let the unique vertex connected to $v$ be called $A$. Note that if every vertex of $G$ were connected to $A$, then the complement of $G$ would be disconnected, as $A$ would have degree 0. Thus, there must exist a vertex $B \in V_G$ s.t. $(A,B) \notin E_G$. Now, examine the edge $(v,B)$. Note that the edge cannot exist since that would violate the condition that $deg(v)=1$. Thus, the edge must either be blocked by triangle freeness or 3-colorability. 
    
    Note that the edge cannot be blocked by triangle freeness as it would require both the edges $(v,C)$ and $(C,B)$ to exist. But for $(v,C)$ to exist, $C = A$ but then we know that the edge $(A,B)$ does not exist. 
    
    Thus, the edge $(v,B)$ must be blocked by 3-colorability. Thus, for any coloring the vertices $v$ and $B$ must be the same color. Without loss of generality, consider vertex $v$ to be red and $A$ to be blue. However, note that since $v$ is not connected to any green vertices, we could recolor $v$ to be green. Thus, there exists a coloring in which $v$ and $B$ are different colors and thus the edge $(v,B)$ cannot be blocked by 3-colorability.

    Thus, the edge $(v,B)$ must exist which means that vertex $v$ cannot have degree 1.
\end{proof}

\begin{mandef}{4} \label{SetDef}
    Let $\Gamma$ be a minimal prime graph complement with vertex $X$ of degree 2. Let the two vertices connected to $X$ be called $A$ and $B$. Let $S_A$ be the set of vertices connected to $A$ but not B; $S_B$ be the set of vertices connected to $B$ but not $A$; $S_Y$ be the set of vertices connected to both $A$ and $B$; and $S_Z$ be the set of vertices connected to neither $A$ or $B$.
\end{mandef}

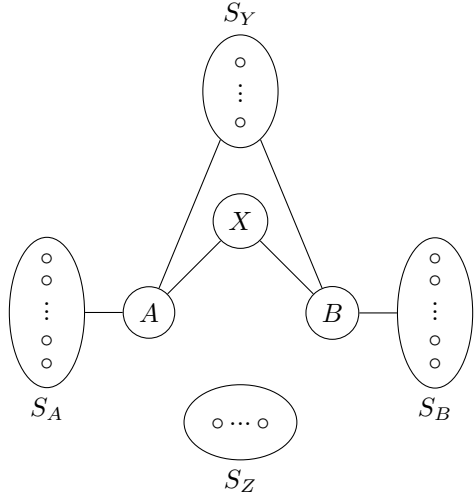
\begin{center}
\begin{tikzpicture}[node distance = 2cm]
\node[circle, draw, fill=white] (1) {$X$};
\node[circle, draw, fill=white] [below left=1cm of 1] (2) {$A$};
\node[circle, draw, fill=white] [below right=1cm of 1] (3) {$B$};
\node[ellipse, minimum height = 1.5cm, minimum width = 1cm, draw, fill=white, label=above:$S_Y$] [above=0.6cm of 1] (4) {};
\node[ellipse, minimum height = 1cm, minimum width = 1.5cm, draw, fill=white, label=below:$S_Z$] [below=1.8cm of 1] (15) {};

\node[ellipse, minimum height = 2cm, minimum width = 1cm, draw, fill=white, label=below:$S_A$] [left=0.5cm of 2] (5) {};

\draw (5 |- 0,-0.5) node (7) {$\circ$};
\draw (5 |- 0,-0.8) node (8) {$\circ$};
\draw (5 |- 0,-1.1) node {.};
\draw (5 |- 0,-1.2) node {.};
\draw (5 |- 0,-1.3) node {.};
\draw (5 |- 0,-1.6) node (9) {$\circ$};
\draw (5 |- 0,-1.9) node (10) {$\circ$};

\node[ellipse, minimum height = 2cm, minimum width = 1cm, draw, fill=white, label=below:$S_B$] [right=0.5cm of 3] (6) {};

\draw (6 |- 0,-0.5) node (11) {$\circ$};
\draw (6 |- 0,-0.8) node (12) {$\circ$};
\draw (6 |- 0,-1.1) node {.};
\draw (6 |- 0,-1.2) node {.};
\draw (6 |- 0,-1.3) node {.};
\draw (6 |- 0,-1.6) node (13) {$\circ$};
\draw (6 |- 0,-1.9) node (14) {$\circ$};

\draw(1 |- 0,1.3) node (16) {$\circ$};
\draw(1 |- 0,1.6) node (17) {.};
\draw(1 |- 0,1.7) node (18) {.};
\draw(1 |- 0,1.8) node (19) {.};
\draw(1 |- 0,2.1) node (21) {$\circ$};

\draw(1 |- 0,-2.7) node (22) {.};
\draw(22 -| 0.1,0) node (23) {.};
\draw(22 -| -0.1,0) node (24) {.};
\draw(24 -| -0.3,0) node (25) {$\circ$};
\draw(23 -| 0.3,0) node (26) {$\circ$};

\draw[-] (1) edge node {} (2);
\draw[-] (1) edge node {} (3);
\draw[-] (2) edge node {} (4);
\draw[-] (2) edge node {} (5);
\draw[-] (3) edge node {} (4);
\draw[-] (3) edge node {} (6);
\end{tikzpicture}
\end{center}


The following lemma shows that the sets $S_A, \phantom{"} S_B,$ and $S_Y$ have no internal edges.
\begin{manlem}{5} \label{intedgesA}
    Let $S \in \{S_A,S_B,S_Y\}$. If $a,b \in S$, then $(a,b) \notin E_{\Gamma}$.
\end{manlem}
\begin{proof}
    We consider the case where $S=S_A$. If the edge $(a,b)$ were to exist, then since the edges $(a,A)$ and $(b,A)$ must exist by Definition \ref{SetDef} there is a triangle $(a,b,A)$. Thus, $(a,b)$ cannot exist.

    Note that the statement is also true by the same logic for two vertices in the sets $S_{B}$ or $S_{Y}$. 
    

\end{proof}

\begin{manlem} {6} \label{Ztri}
    If $v \in S_{Z}$ then the graph $\overline{\Gamma} =(V_{\Gamma},E_{\Gamma} \cup \{(v,X)\})$ is triangle free.
\end{manlem}
\begin{proof}
    Let $a \in V_{\Gamma}$ such that $(v,a)\in E_{\Gamma}$. Note that by the definition of $S_Z$, $v$ cannot be $A$ or $B$. Consequently, the edge $(a,X)$ existing would violate $deg(X) = 2$. Thus, adding the edge $(v,X)$ will not violate triangle-freeness.
\end{proof}
\begin{manlem}{7} \label{intedgesZ}
    If $a,b \in S_{Z}$, then $(a,b) \notin E_{\Gamma}$.
\end{manlem}
\begin{proof}
    Consider the edge $(a,X)$. Note that since $deg(X) = 2$, this edge cannot exist, since $a \ne A$ or $B$. Because minimal prime graph complements are maximal, adding this edge necessarily violates one of their defining conditions. This violated condition can either be triangle freeness or 3-colorability. We know that adding the edge cannot violate triangle freeness by Lemma \ref{Ztri}.

    Because we know that adding the edge $(a,X)$ does not violate triangle freeness, it must violate 3-colorability. That is, $a$ and $X$ must be the same color in every coloring of $\Gamma$; let us call this color red. Assume for the sake of contradiction that there exists $b \in S_{Z}$ such that $(a,b)\in E_{\Gamma}$. We know that because $a$ is red, $b$ cannot be red. Thus, we know that the edge $(b,X)$ is not blocked by 3-colorability since $X$ and $b$ are different colors. Furthermore, the edge cannot violate triangle-freeness by Lemma \ref{Ztri}. Thus, by the maximality condition of a minimal prime graph complement, the edge $(b,X)$ must exist. At the same time, we know that it cannot exist since $deg(X) = 2$, which is a contradiction.
    

\end{proof}

\begin{manlem}{8} \label{YconA}
    If $a \in S_{A}$ and $y \in S_{Y}$ then $(a,y) \notin E_{\Gamma}$.
\end{manlem}
\begin{proof}
    Let $a \in S_{A}$ and $y \in S_{Y}$. Note that by the definitions of $S_A$ and $S_Y$, both $(a,A)$ and $(y,A)$ are edges of $\Gamma$. Thus, if the edge $(a,y)$ were to be in $\Gamma$, then there would be a triangle between vertices $a, y,$ and $A$, which would violate triangle-freeness.
\end{proof}

\begin{manlem}{9} \label{validColor}
     Color the graph in the following way. Color $X$ red; color $A$ green; color all elements of $S_A$ red; color all elements of $S_B$ green; color all elements of $S_Y$ red; color all elements of $S_Z$ blue. This is a valid three coloring of $\Gamma$.
\end{manlem}

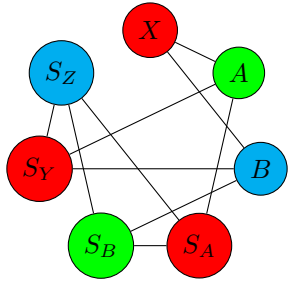
\begin{center}
\begin{tikzpicture}

\node[regular polygon, regular polygon sides=7, minimum size=3cm] (a) {};
    \node[circle, draw, fill=red] (1) at (a.corner 1) {$X$};
    \node[circle, draw, fill=green] (2) at (a.corner 7) {$A$};
    \node[circle, draw, fill=cyan] (3) at (a.corner 6) {$B$};
    \node[circle, draw, fill=red] (4) at (a.corner 5) {$S_A$};
    \node[circle, draw, fill=green] (5) at (a.corner 4) {$S_B$};
    \node[circle, draw, fill=red] (6) at (a.corner 3) {$S_Y$};
    \node[circle, draw, fill=cyan] (7) at (a.corner 2) {$S_Z$};

    \draw[-] (1) edge node {} (2);
    \draw[-] (1) edge node {} (3);
    \draw[-] (2) edge node {} (4);
    \draw[-] (3) edge node {} (5);
    \draw[-] (4) edge node {} (5);
    \draw[-] (4) edge node {} (7);
    \draw[-] (5) edge node {} (7);
    \draw[-] (6) edge node {} (2);
    \draw[-] (6) edge node {} (3);
    \draw[-] (6) edge node {} (7);
\end{tikzpicture}
\end{center}

\begin{proof}
    We will prove this statement in 3 steps; showing that no red vertex is connected to any other; no green vertex is connected to any other green vertex; and no blue vertex is connected to any other blue vertex.
    \begin{enumerate}
        \item RED: Note that since $deg(X) = 2$, no element of $S_A$ or $S_Y$ is connected to $X$. In addition, no element of $S_Y$ can be connected to any element of $S_A$ by Lemma \ref{YconA}. Since by Lemma \ref{intedgesA} there do not exist internal edges within any of the sets $S_A, S_Y,$ or $\{X\}$, no red vertex is connected to any other red vertex.
        \item GREEN: Note that by definition \ref{SetDef}, no element of $S_Z$ is connected to $B$. So, by Lemma \ref{intedgesZ}, there cannot exist any internal edges within the set $S_Z$, there exists no two blue vertices connected to each other.
        \item BLUE: Note that by definition \ref{SetDef} no element of $S_B$ is connected to $A$. So, since $S_B$ has no internal edges by Lemma \ref{intedgesA}, there cannot exist an edge connecting two blue vertices.
    \end{enumerate}

\par
    Thus, there are not 2 vertices in $\Gamma$ of the same color that are connected, and thus this coloring of $\Gamma$ must be a valid 3-coloring.
\end{proof}
\begin{manlem}{10} \label{noZ}
    The set $S_Z$ is empty.
\end{manlem}
\begin{proof}
    Assume for the sake of contradiction that $S_Z$ is non-empty. Let $v\in S_Z$. Now examine the edge $(v, X)$. We know that by Lemma \ref{Ztri}, adding this edge to $\Gamma$ would not violate triangle freeness. We also know that Lemma \ref{validColor} guarantees the existence of a valid coloring of $\Gamma$ such that $v$ and $X$ are different colors. Thus, $(v,X)$ cannot be blocked by 3-colorability. Thus, by the maximality condition, the edge must exist. But since we know that, since $deg(X) = 2$, the edge cannot exist, we have a contradiction. Thus, our assumption that $S_Z$ is nonempty must be false and thus $S_Z$ must be empty.
\end{proof}

\begin{manlem}{11}
    For all $a \in S_A$ and $b \in S_B$, $(a,b) \in E_{\Gamma}$.
\end{manlem}

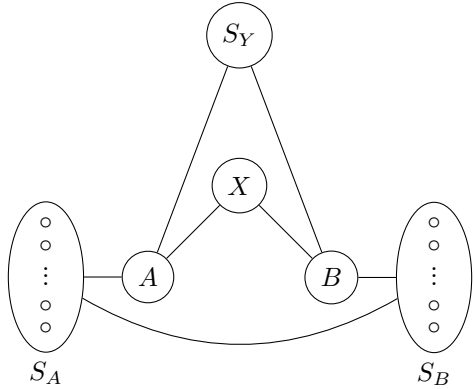
\begin{center}
\begin{tikzpicture}[node distance = 2cm]
\node[circle, draw, fill=white] (1) {$X$};
\node[circle, draw, fill=white] [below left=1cm of 1] (2) {$A$};
\node[circle, draw, fill=white] [below right=1cm of 1] (3) {$B$};
\node[circle, draw, fill=white] [above of=1] (4) {$S_Y$};
\node[ellipse, minimum height = 2cm, minimum width = 1cm, draw, fill=white, label=below:$S_A$] [left=0.5cm of 2] (5) {};

\draw (5 |- 0,-0.5) node (7) {$\circ$};
\draw (5 |- 0,-0.8) node (8) {$\circ$};
\draw (5 |- 0,-1.1) node {.};
\draw (5 |- 0,-1.2) node {.};
\draw (5 |- 0,-1.3) node {.};
\draw (5 |- 0,-1.6) node (9) {$\circ$};
\draw (5 |- 0,-1.9) node (10) {$\circ$};

\node[ellipse, minimum height = 2cm, minimum width = 1cm, draw, fill=white, label=below:$S_B$] [right=0.5cm of 3] (6) {};

\draw (6 |- 0,-0.5) node (11) {$\circ$};
\draw (6 |- 0,-0.8) node (12) {$\circ$};
\draw (6 |- 0,-1.1) node {.};
\draw (6 |- 0,-1.2) node {.};
\draw (6 |- 0,-1.3) node {.};
\draw (6 |- 0,-1.6) node (13) {$\circ$};
\draw (6 |- 0,-1.9) node (14) {$\circ$};

\draw[-] (1) edge node {} (2);
\draw[-] (1) edge node {} (3);
\draw[-] (2) edge node {} (5);
\draw[-] (2) edge node {} (4);
\draw[-] (3) edge node {} (4);
\draw[-] (3) edge node {} (6);
\draw[-] (5) edge [bend right] node[left] {} (6);
\end{tikzpicture}
\end{center}

\begin{proof}
    Let $a \in S_A$ and $b \in S_B$. Assume for the sake of contradiction that the edge $(a,b)$ does not exist. Thus, by the maximality condition, the edge must be blocked by either triangle freeness or 3-colorability. Note that since Lemma \ref{validColor} guarantees the existence of a coloring in which $a$ and $b$ are different colors, the edge cannot be blocked by 3-colorability. Thus, $(a,b)$ must be blocked by triangle freeness. Let $v\in V_{\Gamma}$ be a vertex that blocks $(a,b)$ which means that $(a,v),(v,b) \in E_{\Gamma}$. Since $(a,v)$ and $(a,A) \in E_{\Gamma}$, the edge $(v,A)$ cannot exist. By similar reasoning, the edge $(v,B) \notin E_{\Gamma}$ either. Thus, since $v$ is not connected to either $A$ or $B$, $v \in S_Z$ by definition. Thus, since $S_Z$ is empty by Lemma \ref{noZ}, $v$ cannot exist. Thus, the edge $(a,b)$ cannot be blocked by triangle freeness and thus must exist.
\end{proof}

\begin{manlem} {12}
    Both $S_A$ and $S_B$ are nonempty.
\end{manlem}
\begin{proof}
    First, note that if both $S_A$ and $S_B$ are empty, then all remaining vertices are either $X$, $A$, $B$, or an element of $S_Y$. Note that we can color $A$ and $B$ red since the edge $(A,B)$ cannot exist as it would form a triangle with vertex $X$. In addition, we can color all vertices in the set $S_Y$ and $X$ blue as no $S_Y$ is connected to $X$ $deg(X) = 2$ and there are no internal edges in the set $S_Y$ by Lemma \ref{intedgesA}. Thus, the graph is bipartite and thus, by \cite{gruber} 
    , cannot be the complement of a minimal prime graph.

    Thus, we now handle the case where exactly one of $S_A$ and $S_B$ are empty. Without loss of generality, let that set be $S_A$. Thus, let $b \in S_B$. Note that since $deg(X) = 2$, $b$ cannot be connected to $X$. Also, by Lemma \ref{YconA}, $b$ cannot be connected to any element of $S_Y$. Finally, by Lemma \ref{intedgesA}, $b$ cannot be connected to any other element of $S_B$. Thus, if $b$ were not connected to $A$ then $deg(b) = 1$ which is impossible by Lemma \ref{noDeg1}. Thus, the edge $(b,A)$ must exist. But since we know that it cannot exist by the definition of $S_B$, we have a contradiction. 
\end{proof}

\begin{manlem} {13}
    $\Gamma$ is a reseminant graph.
\end{manlem}
\begin{proof}
    Note that for each $a \in S_A$, $a$ is exactly connected to $A$ and all elements of $S_B$. Thus, all elements of $S_A$ have the exact same adjacencies and thus are vertex duplications of each other. Similarly, all $b \in S_B$ are duplications of each other, and all $y \in S_Y$ are vertex duplications of $X$.
    Removing all duplicated vertices yields a minimal prime graph complement with 5 vertices $X,A,B,a,b$ where $\{(X,A),(A,a),(a,b),(b,B),(B,X)\} = E_{\overline{\Gamma}}$. Note that this is a 5-cycle. Thus, since all vertices in $\Gamma$ are duplicates of some vertex in a 5-cycle, $\Gamma$ is a reseminant graph.
\end{proof}


\section*{Acknowledgments} The authors are grateful to Mathworks 
at Texas State University for supporting this research.

\end{document}